\documentclass[10pt,reqno]{amsart}
\usepackage{amsmath,amsfonts, amssymb, amsthm}
  \usepackage{paralist}
  \usepackage{graphics} 
  \usepackage{epsfig} 
\usepackage{graphicx}  \usepackage{epstopdf}
 \usepackage[colorlinks=true]{hyperref}
\hypersetup{urlcolor=blue, citecolor=red}

\pagespan{559}{566}  
\PII{eISSN: 2688-1594 AIMS (2020)}

\copyrightinfo{2020}{American Institute of Mathematical Sciences}

\newtheorem{theorem}{Theorem}[section]

\newtheorem{lemma}[theorem]{Lemma}

\theoremstyle{definition}

\newtheorem{remark}{Remark}[section]

\def\pmod #1{\ ({\rm{mod}}\ #1)}
\def\Z{\Bbb Z}
\def\N{\Bbb N}

\def\l{\left}
\def\r{\right}
\def\bg{\bigg}
\def\({\bg(}
\def\){\bg)}
\def\t{\text}
\def\f{\frac}

\def\ls{\leq}
\def\gs{\geq}

\def\sm{\setminus}
\def\bi{\binom}

\def\eq{\equiv}

\def\da{\delta}

\def\Proof{\noindent{\it Proof}}

\begin{document}
\hbox{Electron. Res. Arch. 28 (2020), no.\,1, 559-566.}
\medskip

\title[On sums of four pentagonal numbers with coefficients]
      {On sums of four pentagonal numbers\\ with coefficients}

\author[Dmitry Krachun and Zhi-Wei Sun]{Dmitry Krachun and Zhi-Wei Sun$^*$}



\date{January, 2020}  

\subjclass[2010]{Primary: 11B13, 11E25; Secondary: 11D85, 11E20, 11P70.}
 \keywords{Pentagonal numbers, additive bases, ternary quadratic forms.}

\address{
Dmitry Krachun,
St. Petersburg Department of Steklov Mathematical Institute of Russian Academy of Sciences, Fontanka 27, 191023, St. Petersburg, Russia}
\email{dmitrykrachun@gmail.com}
\address{
Zhi-Wei Sun,
Department of Mathematics, Nanjing University, Nanjing 210093,  China}
\email{zwsun@nju.edu.cn}
\thanks{The work is supported by the NSFC(Natural Science Foundation of China)-RFBR(Russian Foundation for Basic Research) Cooperation and Exchange Program (grants NSFC 11811530072 and
RFBR 18-51-53020-GFEN-a). The second author is also supported by the Natural Science Foundation of China (grant no. 11971222).}

\thanks{$^*$ Corresponding author: Zhi-Wei Sun}

\maketitle
\vspace*{-20pt}

%
%
%
%
%

\begin{abstract}
The pentagonal numbers are the integers given by
$p_5(n)=n(3n-1)/2\ (n=0,1,2,\ldots)$.
Let $(b,c,d)$ be one of the triples $(1,1,2),(1,2,3),(1,2,6)$ and $(2,3,4)$.
We show that each $n=0,1,2,\ldots$ can be written as $w+bx+cy+dz$ with $w,x,y,z$ pentagonal numbers,
which was first conjectured by Z.-W. Sun in 2016. In particular, any nonnegative integer
is a sum of five pentagonal numbers two of which are equal; this refines a classical result
of Cauchy claimed by Fermat.
\end{abstract}

\section{Introduction}
\setcounter{theorem}{0}
\setcounter{corollary}{0}
\setcounter{remark}{0}
\setcounter{equation}{0}

For each $m=3,4,5,\ldots$, the {\it polygonal numbers of order $m$} are given by
$$p_m(n)=(m-2)\bi n2+n\ \ (n\in\N=\{0,1,2,\ldots\}).$$
In particular, those $p_5(n)$ with $n\in\N$ are called {\it pentagonal numbers}.
A famous claim of Fermat states that each $n\in\N$  can be written as a sum of $m$ polygonal numbers of order $m$. This was proved by Lagrange for $m=4$ in 1770, by Gauss for $m=3$ in 1796, and
by Cauchy for $m\gs 5$ in 1813. For Cauchy's polygonal number theorem, one may consult Nathanson \cite{N87} and \cite[Chapter 1, pp.\,3-34]{N96} for details. In 1830 Legendre refined Cauchy's polygonal number theorem by showing that for any $m=5,6,\ldots$
every sufficiently large integer is a sum of five polygonal numbers of order $m$ one of which is $0$ or $1$ (cf. \cite[p.\,33]{N96}).

In 2016 Sun \cite[Conjecture 5.2(ii)]{S16} conjectured that each $n\in\N$ can be written as
$$p_5(w)+bp_5(x)+cp_5(y)+dp_5(z)\ \ \t{with}\ w,x,y,z\in\N,$$
provided that $(b,c,d)$ is among the following 15 triples:
\begin{gather*}(1,1,2),(1,2,2),(1,2,3),(1,2,4),(1,2,5),(1,2,6),(1,3,6),
\\(2,2,4),(2,2,6),(2,3,4),(2,3,5),(2,3,7),(2,4,6),(2,4,7),(2,4,8).
\end{gather*}
In 2017, Meng and Sun \cite{MS} confirmed this for $(b,c,d)=(1,2,2),(1,2,4)$.
In this paper we prove the conjecture for
$$(b,c,d)=(1,1,2),\,(1,2,3),\,(1,2,6),\,(2,3,4).$$

\begin{theorem}\label{Th1.1} Each $n\in\N$ can be written as a sum of five pentagonal numbers two of which are equal, that is, there are $x,y,z,w\in\N$ such that
$$n=p_5(x)+p_5(y)+p_5(z)+2p_5(w).$$
\end{theorem}

\begin{remark}\label{Rem1.1} Clearly, Theorem \ref{Th1.1} is stronger than the classical result that any nonnegative integer
is a sum of five pentagonal numbers. In Feb. 2019 the second author even conjectured that
any integer $n>33066$ is a sum of three pentagonal numbers.
\end{remark}

\begin{theorem}\label{Th1.2} Any $n\in\N$ can be written as
$p_5(w)+2p_5(x)+3p_5(y)+4p_5(z)$
with $w,x,y,z\in\N$.
\end{theorem}

\begin{theorem}\label{Th1.3} Let $\da\in\{1,2\}$. Then any $n\in\N$ can be written as
$p_5(w)+p_5(x)+2p_5(y)+3\da p_5(z)$
with $w,x,y,z\in\N$.
\end{theorem}

We will prove Theorems 1.1-1.3 in Sections 2-4 respectively. Our proofs use some known results
on ternary quadratic forms.

Those $p_5(x)=x(3x-1)/2$ with $x\in\Z$ are called {\it generalized pentagonal numbers}.
Clearly,
$$\{p_5(x):\ x\in\Z\}=\l\{\f{n(3n-1)}2:\ n\in\N\r\}\bigcup\l\{\f{n(3n+1)}2:\ n\in\N\r\}.$$
Recently, Ju \cite{Ju}
showed that for any positive integers $a_1,\ldots,a_k$ the set
$$\{a_1p_5(x_1)+\ldots+a_kp_k(x_k):\ x_1,\ldots,x_k\in\Z\}$$
contains all nonnegative integers whenever it contains the twelve numbers
$$1,\ 3,\ 8,\ 9,\ 11,\ 18,\ 19,\ 25,\ 27,\ 43,\ 98,\ 109.$$

The generalized octagonal numbers are those $p_8(x)=x(3x-2)$ with $x\in\Z$. In 2016, Sun \cite{S16} proved that any positive integer can be written as a sum of four generalized octagonal numbers one of which is odd. See also Sun \cite{S20} and \cite{S18} for representations of nonnegative integers in the form
$x(ax+b)/2+y(cy+d)/2+z(ez+f)/2$ with $x,y,z$ integers or nonnegative integers, where $a,c,e$
are positive integers and $b,d,f$ are integers with $a+b,c+d,e+f$ all even.

\section{Proof of Theorem 1.1}
\setcounter{theorem}{0}
\setcounter{corollary}{0}
\setcounter{remark}{0}
\setcounter{equation}{0}

\begin{lemma}\label{Lem2.1} Any positive even number $n$ not in the set $\{5^{2k+1}m:\ k,m\in\N\ \t{and}\ m\eq\pm2\pmod5\}$
can be written as $x^2+y^2+z^2+(x+y+z)^2/2$ with $x,y,z\in\Z$.
\end{lemma}
\Proof. By Dickson \cite[pp.\,112-113]{D39},
\begin{align*}&\N\sm\{x^2+2y^2+10z^2:\ x,y,z\in\Z\}
\\=&\{8m+7:\ m\in\N\}\cup\{5^{2k+1}l:\ k,l\in\N\ \t{and}\ l\eq\pm1\pmod5\}.
\end{align*}
Thus $8n=s^2+2t^2+10z^2$ for some $s,t,z\in\Z$. Clearly, $2\mid s$ and $t\eq z\pmod 2$.
Without loss of generality, we may assume that $t\not\eq z\pmod4$ if $2\nmid z$.
(If $t\eq z\pmod 4$ with $z$ odd, then $-t\not\eq z\pmod4$.)
Write $s=2r$ and $t=2w+z$ with $r,w\in\Z$. Then $2\nmid w$ if $2\nmid z$. Since
$$0\eq 8n=s^2+2(2w+z)^2+10z^2=(2r)^2+12z^2+8w(w+z)\pmod {16},$$
both $r-z$ and $w(w+z)$ are even. If $2\mid z$ then $2\mid w$. Recall that $2\nmid w$ if $2\nmid z$.
So $w\eq z\eq r\pmod2$. Now, both
$x=(r+w)/2$ and $y=(w-r)/2$ are integers. Observe that
\begin{align*}2n=& r^2+2\l(w+\f z2\r)^2+ 10\l(\f z2\r)^2=r^2+(w+z)^2+w^2+2z^2
\\=&(x-y)^2+(x+y)^2+(x+y+z)^2+2z^2
\\=&2x^2+2y^2+2z^2+(x+y+z)^2
\end{align*}
and hence $n=x^2+y^2+z^2+(x+y+z)^2/2$. This ends the proof. \qed

\begin{lemma}\label{Lem2.2} Let $n\in\N$. Suppose that there are $B\in\N$ and $x,y,z\in\Z$ such that  $3\mid n+B$ and
$$\f23(n+B)+B-5B^2=x^2+y^2+z^2+\f{(x+y+z)^2}2\ls B^2.$$
Then $n=p_5(x_0)+p_5(y_0)+p_5(z_0)+2p_5(w_0)$ for some $x_0,y_0,z_0,w_0\in\N$.
\end{lemma}
\Proof. Clearly, $w=-(x+y+z)/2\in\Z$. As $|x|,|y|,|z|,|w|\ls B$, all the numbers
$$x_0=x+B,\ y_0=y+B,\ z_0=z+B,\ w_0=w+B$$
are nonnegative integers. Observe that
\begin{align*}&p_5(x_0)+p_5(y_0)+p_5(z_0)+2p_5(w_0)
\\=&\f{3(x_0^2+y_0^2+z_0^2+2w_0^2)-(x_0+y_0+z_0+2w_0)}2
\\=&\f{3(5B^2+x^2+y^2+z^2+2w^2)-5B}2=\f{2n+5B-5B}2=n.
\end{align*}
This concludes the proof. \qed

\medskip
\noindent{\it Proof of Theorem 1.1}. We can easily verify the desired result for $n=0,\ldots,8891$.
Below we assume that $n\gs 8892$.
 If
\begin{equation}\label{Bn}\f{\sqrt n}3+\f16\ls B\ls\sqrt{\f{2n}{15}}+\f16,
\end{equation}
then
$$\f23(n+B)+B-5B^2\gs\f{15(B-1/6)^2}3+\f{5B}3-5B^2=\f 5{36}>0$$
and
\begin{align*}\f23(n+B)+B-5B^2\ls&\f23\l(3\l(B-\f16\r)\r)^2+\f{5B}3-5B^2
\\=&B^2-\f13\l(B-\f12\r)\ls B^2.
\end{align*}

\medskip
\noindent{\bf Case 1.} $5\nmid n$.

As
$$n\gs\l\lceil\f{3^2}{(\sqrt{2/15}-1/3)^2}\r\rceil=8892,$$
we have
$$\sqrt{\f{2n}{15}}+\f16-\l(\sqrt{\f{n}3}+\f16\r)=\l(\sqrt{\f 2{15}}-\f13\r)\sqrt n\gs 3$$
and hence there is an integer $B$ satisfying \eqref{Bn} with
$B\eq-n\pmod3$. By the above,
$$0\ls\f23(n+B)+B-5B^2=\f{2n+5B}3-5B^2\ls B^2.$$
As the even number $\f23(n+B)+B-5B^2$ is not divisible by $5$, in light of Lemma \ref{Lem2.1} there are $x,y,z\in\Z$ such that
$$\f23(n+B)+B-5B^2=x^2+y^2+z^2+\f{(x+y+z)^2}2.$$
Now, by applying Lemma \ref{Lem2.2} we find that $n=p_5(x_0)+p_5(y_0)+p_5(z_0)+2p_5(w_0)$
for some $x_0,y_,z_0,w_0\in\N$.

\medskip
\noindent{\bf Case 2.} $n=5q$ for some $q\in\N$.

In this case, we can easily verify the desired result when $8892\ls n\ls 222288$.
Below we assume that
$$n\gs222289=\l\lceil\f{15^2}{(\sqrt{2/15}-1/3)^2}\r\rceil.$$
Choose $\da\in\{0,\pm1\}$ such that $1-q-\da\not\eq0,\pm2\pmod5$.
As
$$\sqrt{\f{2n}{15}}+\f16-\l(\sqrt{\f{n}3}+\f16\r)=\l(\sqrt{\f 2{15}}-\f13\r)\sqrt n\gs 15,$$
there is an integer $B$ satisfying \eqref{Bn} such that
$B\eq-n\pmod3$ and $(B-1)^2\eq\da\pmod 5$. Note that
$$\f23(n+B)+B-5B^2=5\l(\f{2q+B}3-B^2\r)$$
and
$$\f{2q+B}3-B^2\eq-\f{2q+B}2-B^2\eq1-q-(B-1)^2\eq1-q-\da\not\eq0,\pm2\pmod5.$$
Thus, by applying Lemmas \ref{Lem2.1} and \ref{Lem2.2} we get that $n=p_5(x_0)+p_5(y_0)+p_5(z_0)+2p_5(w_0)$
for some $x_0,y_,z_0,w_0\in\N$.

In view of the above, we have completed the proof of Theorem \ref{Th1.1}. \qed

\section{Proof of Theorem 1.2}
\setcounter{theorem}{0}
\setcounter{corollary}{0}
\setcounter{remark}{0}
\setcounter{equation}{0}

\begin{lemma}\label{Lem3.1} Let $q\in\N$ with $q$ odd and not squarefree, or $2\mid q$
and $q\not\in\{4^k(16l+6):\ k,l\in\N\}$. Then there are $x,y,z\in\Z$ such that
$$6q=2x^2+3y^2+4z^2+(2x+3y+4z)^2.$$
\end{lemma}
\Proof. By K. Ono and K. Soundararajan \cite {OS}, and Dickson \cite{D27}, the Ramanujan form $x^2+y^2+10z^2$
represents $q$. Write $q=a^2+b^2+10c^2$ with $a,b,c\in\Z$. Then, for
$$x=a+b+2c,\ y=-b+2c,\ z=-3c,$$
we have
$$2x^2+3y^2+4z^2+(2x+3y+4z)^2=6(a^2+b^2+10c^2)=6q.$$
This concludes the proof. \qed

\begin{lemma}\label{Lem3.2} Let $n\in\N$. Suppose that there are $B\in\N$ and $x,y,z\in\Z$ such that
$$\f{2n+10B}3-10B^2=2x^2+3y^2+4z^2+(2x+3y+4z)^2<(B+1)^2.$$
Then $n=p_5(w_0)+2p_5(x_0)+3p_5(y_0)+4p_5(z_0)$ for some $w_0,x_0,y_0,z_0\in\N$.
\end{lemma}
\Proof. Set $w=-(2x+3y+4z)$. As $|w|,|x|,|y|,|z|\ls B$, all the numbers
$$w_0=w+B,\ x_0=x+B,\ y_0=y+B,\ z_0=z+B$$
are nonnegative integers. Observe that
\begin{align*}&p_5(w_0)+2p_5(x_0)+3p_5(y_0)+4p_5(z_0)
\\=&\f{3(w_0^2+2x_0^2+3y_0^2+4z_0^2)-(w_0+2x_0+3y_0+4z_0)}2
\\=&\f{3(10B^2+w^2+2x^2+3y^2+4z^2)-10B}2=\f{2n+10B-10B}2=n.
\end{align*}
This ends the proof.
\qed

\medskip
\noindent{\it Proof of Theorem 1.2}. We can verify the result for $n=0,\ldots,45325137$
directly via a computer. Below we assume that
$$n\gs 45325138=\l\lceil\f{(81-1/6+1/16)^2}{(\sqrt{1/15}-\sqrt{2/33})^2}\r\rceil.$$
Since
$$\sqrt{\f n{15}}+\f16-\l(\sqrt{\frac{2n}{33}}+\f1{16}\r)\gs81,$$
there is an integer $B$ with
$$\sqrt{\frac{2n}{33}}+\f1{16}\ls B\ls \sqrt{\f n{15}}+\f16$$
such that
$$B\eq-9n^3+12n^2-38n\pmod {81}$$
if $n$ is odd, and $B\eq 3n-1\pmod 8$ and $B\eq3n^2-2n\pmod9$ if $n$ is even.
Note that
$$\f{2n+10B}3-10B^2\gs\f{30(B-1/6)^2+10B}3-10B^2=\f 5{18}>0$$
and
\begin{align*}\f{2n+10B}3-10B^2\ls&\f{33(B-1/16)^2+10B}3-10B^2
\\=&B^2+\f{47}{24}B+\f{11}{256}<(B+1)^2.
\end{align*}

Let $q=(n+5B-15B^2)/9$. When $n$ is odd, we can easily see that $q$ is an odd integer divisible by $9$. When $n$ is even, $q$ is an even integer with $q\eq4\pmod8$, and hence
$q\not=4^k(16l+6)$ for any $k,l\in\N$.
By Lemma \ref{Lem3.1}, we can write $6q=(2n+10B)/3-10B^2$ as $2x^2+3y^2+4z^2+(2x+3y+4z)^2$
with $x,y,z\in\Z$. Applying Lemma \ref{Lem3.2}, we see that
$n=p_5(w_0)+2p_5(x_0)+3p_5(y_0)+4p_5(z_0)$ for some $w_0,x_0,y_0,z_0\in\N$.

The proof of Theorem \ref{Th1.2} is now complete. \qed

\section{Proof of Theorem 1.3}
\setcounter{theorem}{0}
\setcounter{corollary}{0}
\setcounter{remark}{0}
\setcounter{equation}{0}

\begin{lemma}\label{Lem4.1} Let $q\in\N$ be a multiple of $9$ with $7\nmid q$ or
$$q\in\{7r:\ r\in\Z\ \t{and}\ r\eq1,2,4\pmod7\}.$$  Then there are $x,y,z\in\Z$ such that
$$6q=x^2+2y^2+3z^2+(x+2y+3z)^2.$$
\end{lemma}
\Proof. Since $9\mid q$ and
$$q\not\in\{7^{2k+1}l:\ k,l\in\N\ \t{and}\ l\eq3,5,6\pmod 7\},$$
by \cite[Theorem 2]{K95} we can write $q$ as $a^2+b^2+7c^2$ with $a,b,c\in\Z$.
For
$$x=6c,\ y=a-b-c,\ z=b-c,$$
we have
$$x^2+2y^2+3z^2+(x+2y+3z)^2=6(a^2+b^2+7c^2)=6q.$$
This concludes the proof. \qed

\begin{lemma}\label{Lem4.2} Let $n\in\N$ and $\da\in\{1,2\}$. Suppose that there are $B\in\N$ and $x,y,z\in\Z$ such that
$$\f{2n+(3\da+4)B}3-(3\da+4)B^2=x^2+2y^2+3\da z^2+(x+2y+3\da z)^2<(B+1)^2.$$
Then $n=p_5(w_0)+p_5(x_0)+2p_5(y_0)+3\da p_5(z_0)$ for some $w_0,x_0,y_0,z_0\in\N$.
\end{lemma}
\Proof. Set $w=-(x+2y+3\da z)$. As $|w|,|x|,|y|,|z|\ls B$, all the numbers
$$w_0=w+B,\ x_0=x+B,\ y_0=y+B,\ z_0=z+B$$
are nonnegative integers. Observe that
\begin{align*}&p_5(w_0)+p_5(x_0)+2p_5(y_0)+3\da p_5(z_0)
\\=&\f{3(w_0^2+x_0^2+2y_0^2+3\da z_0^2)-(w_0+x_0+2y_0+3\da z_0)}2
\\=&\f{3((3\da+4)B^2+w^2+x^2+2y^2+3\da z^2)-(3\da+4)B}2
\\=&\f{2n+(3\da+4)B-(3\da+4)B}2=n.
\end{align*}
This ends the proof.
\qed

\medskip
\noindent{\it Proof of Theorem 1.3 with $\da=1$}. We can verify the desired result for $n=0,1,\ldots,808834880$
directly via a computer. Below we assume that
$$n\gs 808834881=\l\lceil\f{(7\times81+1/48-1/6)^2}{(\sqrt{2/21}-\sqrt{1/12})^2}\r\rceil.$$
Since
$$\sqrt{\frac{2n}{21}}+\f1{6}-\l(\sqrt{\f n{12}}+\f1{48}\r)\gs7\times81,$$
there is an integer $B$ with
$$ \sqrt{\f n{12}}+\f1{48}\ls B\ls \sqrt{\frac{2n}{21}}+\f1{6} $$
such that $B\eq 18n^3+3n^2-35n\pmod{81}$,
and $3n/7+1-(B+1)^2\eq3,5,6\pmod7$ if $7\mid n$.
Such an integer $B$ exists in view of the Chinese Remainder Theorem and the simple observations
\begin{gather*}0-1^2\eq1-3^2\eq 6-0^2\eq6\pmod 7,
\\2-2^2\eq5-0^2\eq5\pmod7,\ 3-0^2\eq 4-1^2\eq3\pmod 7.
\end{gather*}
Note that
$$\f{2n+7B}3-7B^2\gs\f{21(B-1/6)^2+7B}3-7B^2=\f7{36}>0$$
and
$$\f{2n+7B}3-7B^2\ls\f{24(B-1/48)^2+7B}3-7B^2=B^2+2B+\f 1{288}<(B+1)^2.$$

It is easy to see that
$$q:=\f16\l(\f{2n+7B}3-7B^2\r)$$
is an integer divisible by $9$.
If $n=7n_0$ for some $n_0\in\N$, then
\begin{align*}\f q7=&\f16\l(\f{2n_0+B}3-B^2\r)
\\\eq&-\l(\f{9n_0-6B}3-B^2\r)=(B+1)^2-(3n_0+1)\eq1,2,4\pmod 7.
\end{align*}
By Lemma \ref{Lem4.1}, we can write $6q=(2n+7B)/3-7B^2$ as $x^2+2y^2+3z^2+(x+2y+3z)^2$
with $x,y,z\in\Z$. Applying Lemma \ref{Lem4.2} with $\da=1$, we see that
$n=p_5(w_0)+p_5(x_0)+2p_5(y_0)+3p_5(z_0)$ for some $w_0,x_0,y_0,z_0\in\N$.
This completes the proof. \qed

\begin{lemma}\label{Lem4.3} Let $q\in\N$ with $q\not\eq7\pmod 8$ or
$$q\not\in\{5^{2k+1}l:\ k,l\in\N\ \t{and}\ l\eq \pm1\pmod5\}.$$
 Then there are $x,y,z\in\Z$ such that
$$6q=x^2+2y^2+6z^2+(x+2y+6z)^2.$$
\end{lemma}
\Proof. By Dickson \cite[pp.\,112-113]{D39}, we can write $q$ as $a^2+2b^2+10c^2$ with $a,b,c\in\Z$.
For
$$x=2a-b+3c,\ y=-a-b+3c,\ z=-2c,$$
we have
$$x^2+2y^2+6z^2+(x+2y+6z)^2=6(a^2+2b^2+10c^2)=6q.$$
This ends the proof. \qed

\medskip
\noindent{\it Proof of Theorem 1.3 with $\da=2$}. We can verify the desired result for $n=0,1,\ldots,897099188$
directly via a computer.

Below we assume that
$$n\gs 897099189=\l\lceil\f{(360+1/16-1/6)^2}{(\sqrt{1/15}-\sqrt{2/33})^2}\r\rceil.$$
Since
$$\sqrt{\f n{15}}+\f16-\l(\sqrt{\frac{2n}{33}}+\f1{16}\r)\gs5\times8\times9,$$
there is an integer $B$ with
$$\sqrt{\frac{2n}{33}}+\f1{16}\ls B\ls \sqrt{\f n{15}}+\f16$$
such that $B\eq 3n^2-2n\pmod 9$ and $B\eq n^2-n-1\pmod 8$, and $(B-1)^2\not\eq2n_0\pm1,2n_0-2\pmod 5$
if $n=5n_0$ with $n_0\in\N$.
Then
$$q=\f16\l(\f{2n+10B}3-10B^2\r)=\f{n+5B-15B^2}9\in\Z$$
and $q\not\eq7\pmod 8$. If $n=5n_0$ for some $n_0\in\N$, then
\begin{align*}\f q5=&\f{n_0+B-3B^2}9\eq 3B^2-B-n_0\eq\f{B^2-2B}2-n_0
\\=&\f{(B-1)^2-2n_0-1}2\not\eq0,\pm1\pmod5.
\end{align*}
As in the proof of Theorem \ref{Th1.2}, we also have
$$0<6q=\f{2n+10B}3-10B^2<(B+1)^2.$$
Now applying Lemma \ref{Lem4.3} and Lemma \ref{Lem4.2} with $\da=2$, we obtain that
$n=p_5(w_0)+p_5(x_0)+2p_5(y_0)+6p_5(z_0)$ for some $w_0,x_0,y_0,z_0\in\N$.

The proof of Theorem \ref{Th1.3} with $\da=2$ is now complete. \qed


\end{document}